\documentclass[12pt]{amsart}
\usepackage{geometry,graphicx,amscd,color,amsmath,amsfonts,amssymb,paralist}
\usepackage[initials]{amsrefs}

\begin{document}

\title{$L_{p}[0,1] \setminus \bigcup\limits_{q>p} L_{q}[0,1]$ is spaceable for every $p>0$}

\author{G. Botelho \and V. V. F\'avaro \and D. Pellegrino \and J. B. Seoane-Sep\'ulveda}

\address{Faculdade de Matem\'atica, \newline\indent Universidade Federal de Uberl\^{a}ndia, \newline\indent 38.400-902 Ð Uberl\^{a}ndia, Brazil.}\email{botelho@ufu.br, vvfavaro@gmail.com}

\address{Departamento de Matem\'{a}tica, \newline\indent Universidade Federal da Para\'{\i}ba, \newline\indent 58.051-900 - Jo\~{a}o Pessoa, Brazil.}\email{pellegrino@pq.cnpq.br}

\address{Departamento de An\'{a}lisis Matem\'{a}tico,\newline\indent Facultad de Ciencias Matem\'{a}ticas, \newline\indent Plaza de Ciencias 3, \newline\indent Universidad Complutense de Madrid,\newline\indent Madrid, 28040, Spain.}
\email{jseoane@mat.ucm.es}

\thanks{G. Botelho was supported by CNPq Grant 306981/2008-4, V. V. F\'avaro by FAPEMIG Grant CEX-APQ-00208-09, D. Pellegrino by CNPq Grant 620108/2008-8 and J. B. Seoane-Sep\'{u}lveda by the Spanish Ministry of Science and Innovation, grant MTM2009-07848.}

\begin{abstract}
In this short note we prove the result stated in the title; that is, for every $p>0$ there exists an infinite dimensional closed linear subspace of $L_{p}[0,1]$ every nonzero element of which does not belong to $\bigcup\limits_{q>p} L_{q}[0,1]$. This answers in the positive a question raised in 2010 by R. M. Aron on the spaceability of the above sets (for both, the Banach and quasi-Banach cases). We also complete some recent results from \cite{BDFP} for subsets of sequence spaces.
\end{abstract}


\maketitle

As it has become an standard notion in the last years, given a topological vector space $X$, we say that a subset $M \subset X$ is {\em spaceable} (see  \cite{AGS}) if there exists an infinite dimensional closed linear subspace $Y \subset M \cup \{0\}$. Very recently, it was proved in \cite{BDFP} that, for every $p>0$, the set  $\ell_{p} \setminus \textstyle\bigcup\limits_{0<q<p} \ell_{q}$ is spaceable. As a consequence of a lecture delivered by the second author at an international conference held in Valencia (Spain) in the summer of 2010, R. M. Aron asked the question of whether a similar result to \cite[Corollary 1.7]{BDFP} would hold for $L_p$-spaces. The aim of this note is to answer R. M. Aron's question in the positive by means of a constructive procedure and some classical Real Analysis and Linear Algebra techniques. \\

\noindent \textbf{Theorem.\,} {\em $L_{p}[0,1] \setminus \bigcup\limits_{q>p} L_{q}[0,1]$ is spaceable for every $p>0$.}\\
\noindent {\sc Proof.} Let us first consider the following representation of the semi-open interval $[0,1)$ as a disjoint union of intervals:
$$\lbrack 0,1)=\left[0,1-1/2\right)  \cup\left[1-1/2,1-1/4\right)  \cup\left[1-1/4,1-1/8\right)  \cup \dots= \textstyle\bigcup\limits_{n =1}^{\infty} I_n,$$
where $I_n := [a_n,b_n)=\left[1-\frac{1}{2^{n-1}},1-\frac{1}{2^{n}}\right)$. Notice that, for every $n \in \mathbb{N}$ and every $x \in I_{n}$, there is a
unique $x_{n}\in\left[  0,1\right)$ such that
$$x=(1-x_{n})a_{n}+x_{n}b_{n}.$$
Now, given $p>0$, let us fix a function $f \in L_{p}[0,1]-\bigcup_{q>p}L_{q}[0,1]$, and define a sequence of functions $(f_{n})_{n=1}^{\infty}$, with $f_{n}\colon\lbrack0,1]\longrightarrow\mathbb{R}$, as follows:
$$f_{n}(x)=\left\{
\begin{array}
[c]{cl}
f(x_{n}) & \mathrm{if~}x\in I_{n},\\
0 & \mathrm{if~}x\notin I_{n}.
\end{array}
\right.$$
The geometric idea is to reproduce the graph of $f$ on the interval $I_{n}$.
By construction, we have that $\left\Vert f_{n}\right\Vert _{L_{p}}\leq\left\Vert
f\right\Vert _{L_{p}}$ for every $n\in\mathbb{N}.$ Also, the functions
$f_{n}$ are linearly independent (they have disjoint supports) and $$\mathrm{span}\{f_{n}:n\in\mathbb{N}\}\subset L_{p}[0,1]\setminus \textstyle\bigcup\limits_{q>p}L_{q}[0,1].$$
The latter proves that $L_{p}[0,1] \setminus \bigcup_{q>p}L_{q}[0,1]$ is $\aleph_{0}$-lineable (as it was seen in \cite{MPPS} for the case $p>1$). In order to go further, the strategy shall be to define a bounded linear and injective operator $T\colon F \longrightarrow L_{p}[0,1]$, where $F$ is a Banach space, and such that $\overline{T(F)}\cap L_{q}[0,1]=\{0\}$ for every $q>p$. This shall prove $L_{p}[0,1]-\bigcup_{q>p}L_{q}[0,1]$ is spaceable. Indeed, if $(\alpha_{j})_{j=1}^{\infty}\in\ell_{s},$ where $s=1$ if $p\geq1$ and $s=p$ if $0<p<1$, we have
$$
\sum_{n=1}^{\infty}\left\Vert \alpha_{n}f_{n}\right\Vert _{L_{p}}^{s} =\sum_{n=1}^{\infty}\left\vert \alpha_{n}\right\vert ^{s}\left\Vert
f_{n}\right\Vert _{L_{p}}^{s}\leq\sum_{n=1}^{\infty}\left\vert \alpha_{n}\right\vert ^{s}\left\Vert f\right\Vert _{L_{p}}^{s}=\left\Vert
f\right\Vert _{L_{p}}^{s}\left\Vert (\alpha_{j})_{j=1}^{\infty}\right\Vert_{s}^{s}<+\infty.
$$
Since $L_{p}[0,1]$ is a Banach space for $p>1$ and a quasi-Banach space for $0<p<1$, it follows that $\sum_{n=1}^{\infty}\alpha_{n}f_{n}\in L_{p}[0,1]$ and, thus,
$$T\colon\ell_{s}\longrightarrow L_{p}[0,1]~~,~~T((\alpha_{j})_{j=1}^{\infty})=\sum_{n=1}^{\infty}\alpha_{n}f_{n}$$
is a well defined linear operator. Suppose that $\sum_{n=1}^{\infty}\alpha_{n}f_{n}=0.$ Then (given any $k \in \mathbb{N}$) we have that, for all $x\in I_{k},$
$$\alpha_{k}f_{k}\left(  x\right)  =\sum_{n=1}^{\infty}\alpha_{n}f_{n}\left(x\right) =0$$
and, since each $f_{k}\neq 0$, it follows that $\alpha_{k}=0.$ Hence $T$ is injective and $T\left(  \ell_{s}\right)  $ is a linear subspace of
$L_{p}[0,1].$ Considering the closure $\overline{T\left(  \ell_{s}\right)  }$ of $T\left(  \ell_{s}\right)  $ in $L_{p}[0,1],$ let us prove that
$\overline{T\left(  \ell_{s}\right)  }\cap L_{q}[0,1]=\{0\}$ for every $q>p.$ Indeed, let $g\in\overline{T\left(  \ell_{s}\right)}\setminus \{0\}$. There exist sequences $\left(a_{i}^{(k)}\right)_{i=1}^{\infty}\in\ell_{s}$ ($k\in\mathbb{N}$) such that $g=\lim_{k\rightarrow\infty}T\left(  \left(a_{i}^{(k)}\right)_{i=1}^{\infty}\right)  $ in $L_{p}[0,1].$
Thus, we have
$$
\displaystyle \lim_{k\rightarrow\infty} \int_{0}^{1} \left \vert \sum_{n=1}^{\infty}a_{n}^{(k)}f_{n}\left(x\right) -g(x)\right\vert ^{p}dx=0,
$$
which, in particular, implies that
$$
\displaystyle \lim_{k\rightarrow\infty} \int_{0}^{\frac{1}{2}} \left\vert a_{1}^{(k)}f_{1}\left(  x\right)  -g(x)\right\vert ^{p}dx=0.
$$
Next, and by means of a subsequence (if needed), we obtain that
$$\lim_{k\rightarrow\infty}a_{1}^{(k)}f_{1}\left(  x\right)  = g(x)\text{ a.e. } x\in [0,1/2).$$
Clearly, the set  $S = \left\{  x\in\lbrack0,\frac{1}{2}) : f_{1}\left(  x\right)
\neq0\right\}$ is not of measure zero. Thus, given $x' \in S$, we have
$$\lim_{k\rightarrow\infty}a_{1}^{(k)}=\frac{f_{1}(x')}{g(x')} = \eta \neq 0,$$
concluding that
$$
g(x)=\eta f_{1}(x) \text{ a.e. } x\in [0,1/2).
$$
which implies that $g\notin L_{q}\left[0,1/2\right]$ (regardless of the $q>p$), finishing the proof.\hfill $\blacksquare$\\

\noindent \textbf{Remark.\,} Let us point out that the previous result is the best possible in terms of dimension, since $L_p[0,1]$ is a $\mathfrak{c}$-dimensional linear space (with $\mathfrak{c}$ denoting the continuum). Also, notice that in the above proof, we actually have that $g\notin L_{q}(I_n)$ regardless of the $q>p$ and $n \in \mathbb{N}$.

\end{document}